\documentclass[english,11pt]{article}%
\usepackage{bbm,ulem}
\usepackage{graphicx}
\usepackage{subcaption}
\usepackage{makeidx}
\usepackage{amssymb}
\usepackage[utf8]{inputenc}
\usepackage{amsfonts}
\usepackage{amsmath}
\usepackage{amssymb}
\usepackage{amsthm}
\usepackage{url}
\usepackage[colorlinks,citecolor=blue]{hyperref}
\usepackage{mathabx}
\usepackage{wrapfig}
\usepackage{dsfont}
\usepackage{resizegather}
\usepackage{yfonts}
\usepackage{authblk}
\usepackage{geometry}
\usepackage[dvipsnames]{xcolor}

\DeclareMathOperator*{\argmin}{arg\,min}

\providecommand{\U}[1]{\protect\rule{.1in}{.1in}}

\usepackage{xcolor}

\newtheorem{prop}{Proposition}[section]
\newtheorem{cor}[prop]{Corollary}

\newtheorem{lem}[prop]{Lemma}

\newtheorem{theo}[prop]{Theorem}

\newcommand{\tr}{\mbox{\rm Tr}}

\newcommand{\EE}{\mathbb{E}}

\newcommand{\LL}{\mathbb{L}}

\newcommand{\RR}{\mathbb{R}}

\newcommand{\TT}{\mathbb{T}}
\newcommand{\UU}{\mathbb{U}}
\newcommand{\VV}{\mathbb{V}}

\newcommand{\Da}{ {\cal D }}
\newcommand{\La}{ {\cal L }}

\newcommand{\Ka}{ {\cal K }}

\newcommand{\Ma}{ {\cal M }}

\newcommand{\Ha}{ {\cal H }}
\newcommand{\Ja}{ {\cal J }}
\newcommand{\Pa}{ {\cal P }}

\newcommand{\point}{\mbox{\LARGE .}}
\newcommand{\cqfd}{\hfill\blbx \\}
\def\blbx{\hbox{\vrule height 5pt width 5pt depth 0pt}\medskip}

\def \RR{\mathbb{R}}

\def \EE{\mathbb{E}}

\def \LL{\mathbb{L}}

\newcommand{\cchi}{\protect\raisebox{2pt}{$\chi$}}

\numberwithin{equation}{section}

\newcommand{\vertiii}[1]{{\left\vert\kern-0.25ex\left\vert\kern-0.25ex\left\vert #1
    \right\vert\kern-0.25ex\right\vert\kern-0.25ex\right\vert}}

\usepackage{amsopn}
\begin{document}

  \title{Entropic continuity bounds for conditional covariances
  with applications to Schr\" odinger and Sinkhorn bridges}

\author{P. Del Moral
\thanks{Centre de Recherche Inria Bordeaux Sud-Ouest, Talence, 33405, France. {\footnotesize E-Mail:\,} \texttt{\footnotesize pierre.del-moral@inria.fr}}}

\maketitle
  \begin{abstract}  
  The article presents new entropic continuity bounds for conditional expectations and conditional covariance matrices.
  These bounds are expressed  in terms of the relative entropy between different coupling distributions. 
   Our approach combines Wasserstein coupling with quadratic transportation cost
inequalities.  We illustrate the impact of these results
in the context of entropic optimal transport problems.

The entropic continuity theorem presented in the article allows to estimate 
  the conditional  expectations and the conditional covariances of  Schr\" odinger and Sinkhorn transitions in terms of the relative 
  entropy between the corresponding  bridges.
These entropic continuity bounds turns out to be a very useful tool for obtaining remarkably simple
proofs of the exponential decays 
of the gradient and the Hessian of  Schr\"odinger and Sinkhorn bridge potentials.
\\

\end{abstract}


\section{An entropic continuity theorem}

Conditional expectations and conditional covariance matrices play a crucial role in a 
variety of areas in applied mathematics  including in statistical inference, nonlinear filtering and data assimilation 
as well as in entropic optimal transport. 

Our aim is to estimate the continuity properties for these conditional quantities
 in terms of the relative entropy between different joint distributions. 
 Our approach combines Wasserstein coupling with quadratic transportation cost
inequalities.
  To the best of our knowledge,
 the  relative entropy continuity bounds presented in this article
are the first result of this type for this class of conditional expectations and conditional covariance matrices.

 We illustrate the impact of these continuity bounds
in the context of entropic optimal transport problems.
Entropic optimal transport methodologies, including Schr\"odinger bridges and Sinkhorn algorithm
 have become state-of-the-art tools in generative modeling and machine learning, see for instance~\cite{bortoli-heng,cuturi,kolouri,peyre} and references therein.

 In this context, an important problem is to estimate the conditional mean and the conditional covariances associated with  these
 optimal coupling distributions~\cite{adm-24,dp}. 
 As we shall see in the sequel, these conditional quantities are also directly related to the gradient and the Hessian of 
 the Schr\" odinger and Sinkhorn bridge potentials~\cite{dp,gt},  see also formulae (\ref{lem-u-mc}) 
 (\ref{lem-v-mc}) in the present article. The entropic continuity theorem presented in this article  (Theorem~\ref{theo-1}) provides an estimate of 
  these conditional  quantities in terms of the entropy between Schr\" odinger and Sinkhorn bridges (see Corollary~\ref{cor-even} and Corollary~\ref{cor-odd}).

   The exponential decays of the relative entropy between Schr\" odinger and Sinkhorn bridges when the number of iterations
 tends to $\infty$ is a very active research area in applied probability and machine learning. It is clearly out of the scope of this note to review 
 this subject, we rather refer to the books~\cite{nutz,peyre} as well as the  recent article articles~\cite{adm-24,adm-25,chiarini,durmus,dp,gt} and the references therein.
  
  In the present article, we illustrate these entropic exponential decays in the context of "log-concave 
at infinity" type marginals measures (in the sense that the curvature estimates discussed in (\ref{reg-UV}) hold, see for instance Theorem~\ref{theo-princ}).  

The refined analysis of more general models require to develop more sophisticated 
mathematical techniques. For instance, the article~\cite{adm-25} presents a novel
approach to Sinkhorn exponential convergence based on Lyapunov techniques and
contraction coefficients on weighted Banach spaces. This analysis applies
to  Gaussian models and statistical finite mixture models.
 We also refer the reader to the stochastic
control representation of entropic plans and back propagation techniques along
Hamilton Jacobi equations used in~\cite{chiarini,conforti-ptrf,gt} to handle general classes of
 weakly semi-concave models including situations where 
 the gradient/Hessian of the potentials 
 are not defined.

To state with some precision our main result we need to introduce some terminology.
 For a given Markov transition $\Ka(x,dy)$ on $\RR^d$ we denote by 
 $ m_{\Ka}(x)$ and $ \Sigma_{\Ka}(x)$ the conditional expectation and the conditional
 covariance matrix defined for any $x\in\RR^d$ by
 $$
 m_{\Ka}(x):=\int~\Ka(x,dy)~y\quad\mbox{\rm and}\quad
 \sigma_{\Ka}(x):=\int~\Ka(x,dy)~(y- m_{\Ka}(x))(y- m_{\Ka}(x))^{\prime}
 $$
 In the above display, points $x,y\in\RR^d$ are represented by column vectors,
 integrals of vector and matrix valued functions are defined by 
 vector and matrices with entries given by 
   integrals w.r.t. the corresponding coordinate.

For a given probability measure $\mu$ and a function $f$ on $\RR^d$, 
 we also consider the product measure
\begin{eqnarray}
 \Ka(f)(x)&:=&\int~\Ka(x,dy)~f(y)\nonumber\\
 (\mu\Ka)(dy)&:=&\int~\mu(dx)\, \Ka(x,dy)
\quad \mbox{\rm and}\quad
 (\mu\times \Ka)(d(x,y)):=\mu(dx)\, \Ka(x,dy)
\label{prod-nott}
\end{eqnarray}
We denote by $\Da_{2}(\mu,\eta)$ and $ \Ha(\mu~|~\eta)$ the $2$-Wasserstein distance between probability measures $\mu$ and $\eta$ 
and the relative entropy of $\mu$ with respect to some $\eta$. 

We are now in position to state the main result of the article.
 \begin{theo}\label{theo-1}
Let $\Ka$ be a Markov transition satisfying the   quadratic transportation cost inequality $\TT_2(\rho)$ for some $\rho>0$. In this situation, for any
probability measure $\mu$ and any Markov transition $\La$ we have
\begin{equation}\label{eq-0}
\Da_2\left(\mu \La,\mu \Ka\right)^2\leq \int\mu(dx)~\Da_2\left(\delta_x\La,\delta_x\Ka\right)^2
\leq 2\rho~\Ha\left(\mu\times\La~|~\mu\times\Ka\right)
\end{equation}
In addition, we have the bias and covariance entropic bounds estimates
 \begin{eqnarray}
 \vertiii{m_{\Ka}-m_{\La}}_{2,\mu}^2&\leq &
 2\rho~\Ha\left(\mu\times\La~|~\mu\times\Ka\right)\label{eq-1}\\
 \vertiii{\sigma_{\Ka}-\sigma_{\La}}_{1,\mu}
&\leq& 4\rho~\Ha\left(\mu\times\La~|~\mu\times\Ka\right)+c_{\mu,\La}~\left(8\rho~\Ha\left(\mu\times\La~|~\mu\times\Ka\right)\right)^{1/2}\label{eq-2}
 \end{eqnarray}
 with the parameter
$
c_{\mu,\La}^2:=\int\mu(dx)\tr(\sigma_{\La}(x))
$ defined in terms of the trace $\tr(\sigma_{\La}(x))$ of the conditional covariance matrix $\sigma_{\La}(x)$. 
 \end{theo}

\section{Some basic notation}

The spectral norm of a matrix $v$ is defined by $\Vert v\Vert_2=\sqrt{\ell_{\tiny max}(v^{\prime}v)}$, with  $v^{\prime}$ the transpose.  The Frobenius norm
 $\Vert v\Vert_F=\sqrt{\tr(v^{\prime}v)}$ with the trace $\tr(v)$ of the matrix $v$. We usually represent points $x=(x^j)_{1\leq j\leq d}\in\RR^d$ by $d$-dimensional column vectors and $1 \times d$ matrices. In this notation, the Frobenius and the spectral norm $\Vert x\Vert_F=\Vert x\Vert_2=\Vert x\Vert:=\sqrt{x^{\prime}x}$ coincides with the Euclidean norm.
 
Consider some probability measures $\mu$ and $\eta$ on $\RR^d$ and some parameter $p\geq 1$.  
Let $\Pi(\mu,\eta)$ be the convex subset of probability measures $P(d(x,y))$ on $
(\RR^d\times\RR^d)$ with marginal $\mu(dx)$ w.r.t. the first $x$-coordinate and marginal $\eta(dy)$ w.r.t. the second $y$-coordinate.

 The $p$-th Wasserstein distance between $\mu$ and $\eta$ and the $\LL_p(\mu)$-norm of a measurable map
$x\in\RR^d\mapsto v(x)\in\RR^{d_1\times d_2}$ are defined by
\begin{equation}\label{def-c2}
\Da_{p}(\mu,\eta)^p:= \inf_{P \in \Pi(\mu,\eta)}~\int~P(d(x,y))~ \Vert x-y\Vert^p~
\quad\mbox{\rm and}\quad
 \vertiii{v}^p_{p,\mu}:=\int\mu(dx)~\Vert v(x)\Vert_F^p
\end{equation}
 The relative entropy (a.k.a. Kullback-Leibler divergence and $I$-divergence) and the Fisher information of  $\mu$ with respect to some $\eta\gg \mu$ is defined by
\begin{equation}\label{def-entropy}
 \Ha(\mu~|~\eta):=\mu\left(\log{\left({d\mu}/{d\eta}\right)}\right)
\quad\mbox{\rm and}\quad
\Ja(\mu~|~\eta)=\mu\left(\Vert\nabla \log{(d\mu/d\eta)}\Vert^2\right)
\end{equation}
 We also use the convention $  \Ha(\mu~|~\eta)=\infty=\Ja(\mu~|~\eta)$ when $\mu\not\ll \eta$.  
 
 We say that a measure $\mu\in\Ma_1(\RR^d)$ satisfies the log-Sobolev inequality $LS(\rho)$  if for any $\nu\in\Ma_1(\RR^d)$  we have
\begin{eqnarray}
 \Ha(\nu~|~\mu)
&\leq &\frac{\rho}{2}~
\Ja(\nu~|~\mu)\label{K-log-sob}
\end{eqnarray}
We say that a Markov transition $\Ka$ satisfies  the log-Sobolev inequality $LS(\rho)$ when $\delta_x\Ka$ satisfies a   log-Sobolev inequality $LS(\rho)$ for any $x\in\RR^d$.
 
 We also say that a probability measure $\mu$ on $\RR^d$ satisfies the quadratic transportation cost inequality $\TT_2(\rho)$ with constant $\rho>0$ if for any probability measure $\nu$  we have
\begin{eqnarray}
\frac{1}{2}~\Da_{2}(\nu,\mu)^2&\leq& \rho~ \Ha(\nu~|~\mu)\label{K-talagrand}
\end{eqnarray}
We say that a Markov transition $\Ka$ satisfies the quadratic transportation cost inequality $\TT_2(\rho)$ when $\delta_x\Ka$ satisfies a quadratic transportation cost inequality $\TT_2(\rho)$ for any $x\in\RR^d$.

A theorem by Otto and Villani, Theorem 1 in~\cite{otto-villani}, ensures that the  log-Sobolev inequality $LS(\rho)$
implies the the quadratic transportation cost inequality $\TT_2(\rho)$. Log-Sobolev inequalities are  a well studied technique for the stability analysis of Markov processes. The class of measures satisfying these inequalities is quite large, see for instance~\cite{toulouse-team,bgl,cattiaux-14,cattiaux, bolbeault,ledoux-99,ledoux,royer} 
and the pioneering articles by Gross~\cite{gross,gross-2}.

Next, we illustrate these conditions with the rather well known example of uniformly log-concave at infinity transitions.
Equip $\RR^d$ with the Lebesgue measure $\lambda(dx)$ and
denote by  $K_{W}(x,dy)$  the
Markov transition  
\begin{equation}
K_W(x,dy):=e^{-W(x,y)}~\lambda(dy)\label{ref-KW}
\end{equation}
for some smooth function $W(x,y)$ on $(\RR^d\times\RR^d)$. 
We now recall a simple criterion for $K_{W}$ to satisfy a $LS(\rho)$ inequality. 
Assume there exists some $a>0$, $b\geq 0$ and $\delta\geq 0$ such that for any $x$ and
$y\in\RR^d$ we have
$$
  \nabla_2^2W(x,y)\geq a~1_{\Vert y\Vert \geq  \delta}~I-b 
~1_{\Vert y\Vert< \delta}~I
$$
In the above display, $  \nabla_2^2W(x,y)$ stands for the Hessian of $W(x,y)$ w.r.t. the second $y$-coordinate. In this case, $K_W$ satisfies the $LS(\rho)$ and the $\TT_2(\rho)$-inequalities (\ref{K-log-sob}) with some parameter $\rho=\rho(a,b,\delta)$ that depends on 
$(a,b,\delta)$. In addition, when $\delta=0$ we can choose $\rho=1/a$. For a more detailed discussion on these curvature criteria we refer to Section 4.2 and Section 7.21 in~\cite{dp}

\section{Proof of Theorem~\ref{theo-1}}
 
Consider the cross covariance matrix $C_{X,Y}$ associated with some random vectors $(X,Y)$ taking values in $(\RR^d\times\RR^d)$ defined by
$$
C_{X,Y}:=\EE\left((X-\EE(X))(Y-\EE(Y))^{\prime}\right)
$$
\begin{lem}
For any random variables $(X,Y)$ with marginal distributions $(\nu^X,\nu^Y)$  we have
the bias and the covariance estimates
\begin{eqnarray}
\Vert \EE(X)-\EE(Y)\Vert_2&\leq& \Da_2(\nu^X,\nu^Y)\label{bias-diff}
\\
\Vert C_{X,X}-C_{Y,Y}\Vert_{F}&\leq& 2~\Da_2(\nu^X,\nu^Y)^2+2~\Da_2(\nu^X,\nu^Y)~(\tr(C_{Y,Y}))^{1/2}\label{cov-diff}
\end{eqnarray}
\end{lem}
\proof
Recalling that $\EE(Z)^{\prime}~\EE(Z)\leq \EE(Z^{\prime}Z)=\EE(\Vert Z\Vert^2)$
we have the bias estimate
  \begin{eqnarray*}
\Vert \EE(X)-\EE(Y)\Vert_2^2= (\EE(X-Y))^{\prime}~ \EE(X-Y)\leq \EE(\Vert X-Y\Vert_2^2)
 \end{eqnarray*}
Note that the above estimate is valid for any coupling of the random variables $(X,Y)$.
We end the proof of (\ref{bias-diff}) choosing the $2$-Wasserstein coupling between  the marginal distributions $(\nu^X,\nu^Y)$.

Now we come to the proof of (\ref{cov-diff}). We have
$$
C_{X,X}-C_{Y,Y}=C_{X-Y,X-Y}+C_{X-Y,Y}+C_{Y,X-Y}
$$
Also note that
$$
\begin{array}{l}
C_{X-Y,X-Y}=\EE((X-Y)(X-Y)^{\prime})-(\EE(X)-\EE(Y))
(\EE(X)-\EE(Y))^{\prime}\\
\\
\Longrightarrow \Vert C_{X-Y,X-Y}\Vert_F\leq \Vert \EE((X-Y)(X-Y)^{\prime})\Vert_F+
\Vert (\EE(X)-\EE(Y))
(\EE(X)-\EE(Y))^{\prime}\Vert_F
\end{array}
$$
Recalling that for any positive definite matrix $C$ we have  $\Vert C\Vert^2_{F}=\tr(C^2)\leq (\tr(C))^2$ we check that
$$
\Vert C_{X-Y,X-Y}\Vert_F\leq \EE(\Vert X-Y\Vert_2^2)+
\Vert \EE(X)-\EE(Y)\Vert_2^2\leq 2~\EE(\Vert X-Y\Vert_2^2)
$$
On the other hand, by Cauchy-Schwartz inequality we have
$$
\Vert C_{Y,X-Y}\Vert^2_F=\Vert C_{X-Y,Y}\Vert^2_F=\sum_{i,j}\EE((X_i-Y_i)(Y_j-\EE(Y_j)))^2\leq 
\EE(\Vert X-Y\Vert_2^2)~\tr(C_{Y,Y})
$$
We conclude that
$$
\Vert C_{X,X}-C_{Y,Y}\Vert_F\leq 2~\EE(\Vert X-Y\Vert_2^2)+2
\EE(\Vert X-Y\Vert_2^2)^{1/2}~(\tr(C_{Y,Y}))^{1/2}
$$
Note that the above estimate is valid for any coupling of the random variables $(X,Y)$.
We end the proof of (\ref{cov-diff}) choosing the $2$-Wasserstein coupling between  the marginal distributions $(\nu^X,\nu^Y)$.
This ends the proof of the lemma.\cqfd

Now we come to the proof of Theorem~\ref{theo-1}.

{\bf Proof of Theorem~\ref{theo-1}:} 
We have
$$
\begin{array}{l}
\displaystyle\pi_x:=
\argmin_{P \in \Pi(\delta_x\La,\delta_x\Ka)}~\int~P(d(y_1,y_2))~ \Vert y_1-y_2\Vert^2
\\
\\
\displaystyle\Longleftrightarrow\Da_2\left(\delta_x\La,\delta_x\Ka\right)^2=\int~\pi_x(d(y_1,y_2))~ \Vert y_1-y_2\Vert^2
\end{array}
$$
This implies that
$$
\begin{array}{l}
\displaystyle\pi_{\mu}(d(y_1,y_2)):=\int \mu(dx)~\pi_x(d(y_1,y_2))\in \Pi(\mu\La,\mu\Ka)\\
\\
\displaystyle\Longrightarrow
\Da_2\left(\mu \La,\mu \Ka\right)^2\leq \int\pi_{\mu}(d(y_1,y_2))~\Vert y_1-y_2\Vert^2=
\int\mu(dx)~\Da_2\left(\delta_x\La,\delta_x\Ka\right)^2
\end{array}$$
Since $\Ka$ satisfies the $\TT_2(\rho)$ inequality we have
$$
 \Da_2(\delta_x\Ka,\delta_x\La)^2\leq 2\rho~\Ha\left(\delta_x\La~|~\delta_x\Ka\right)
 $$
 We end the proof of (\ref{eq-0}) using the fact that
$$
\int \mu(dx)~\Ha\left(\delta_x\La~|~\delta_x\Ka\right)=
\Ha\left(\mu\times\La~|~\mu\times\Ka\right)
$$
By (\ref{bias-diff}) we have
 $$
\displaystyle\vertiii{m_{\Ka}-m_{\La}}_{2,\mu}^2\leq \int\mu(dx)~\Da_2(\delta_x\Ka,\delta_x\La)^2
$$
 In the same vein, combining (\ref{cov-diff}) with Cauchy-Schwartz inequality we have
 $$
 \begin{array}{l}
\displaystyle2^{-1}
\vertiii{\sigma_{\Ka}-\sigma_{\La}}_{1,\mu}
\displaystyle\leq \int\mu(dx)~\Da_2(\delta_x\Ka,\delta_x\La)^2+c_{\mu,\La}~\left(\int\mu(dx)~\Da_2(\delta_x\Ka,\delta_x\La)^2\right)^{1/2}
\end{array}
$$
The end of the proofs of (\ref{eq-1})  and (\ref{eq-2}) now follows word-for-word the same lines of arguments as above, thus it is skipped. This ends the proof of the theorem.\cqfd

\section{Entropic optimal transport}
\subsection{Schr\" odinger  and Sinkhorn bridges}
The Schr\" odinger bridge from $\mu$ to 
$\eta$ with respect to a given reference probability distribution $ P$ on $(\RR^d\times\RR^d)$ is defined by
\begin{equation}\label{ob-eq}
P_{\mu,\eta}:= \argmin_{Q\,\in\, \Pi(\mu,\eta)}\Ha(Q~|~P)
 \end{equation}
It is implicitly  assumed there exists some $Q\in \Pi(\mu,\eta)$ such that 
$\Ha(Q~|~P)<\infty$.
This condition ensures  the existence of a Schr\" odinger bridge distribution $P_{\mu,\eta}$
 solving (\ref{ob-eq})  (cf. the seminal article by Csisz\'ar~\cite{csiszar-2}, as well as Section 6 in the Lecture Notes by Nutz~\cite{nutz}, see also the survey article by L\'eonard~\cite{leonard} and numerous references therein). 
 
  In contrast with conventional optimal transport maps with limited applicability in high dimensions,  Schr\"odinger bridges (\ref{ob-eq})
 can be solved using the celebrated Sinkhorn's algorithm, a.k.a. iterative proportional fitting procedure~\cite{cuturi,genevay-cuturi-2,genevay,sinkhorn,sinkhorn-2,sinkhorn-3}.   Next we briefly outline the basic principles of the algorithm.
Sinkhorn algorithm starts from the reference measure $\Pa_0=P$ and solves sequentially the following entropic transport problems
\begin{equation}\label{sinhorn-entropy-form}
\Pa_{2n+1}:= \argmin_{\Pa\in \Pi(\point,\eta)}\Ha(\Pa~|~\Pa_{2n})\quad \mbox{\rm and}\quad
\Pa_{2(n+1)}:= \argmin_{\Pa\in  \Pi(\mu,\point)}\Ha(\Pa~|~\Pa_{2n+1}).
\end{equation}
where $\Pi(\point,\eta)$, respectively $ \Pi(\mu,\point)$, stands for the distributions  $Q$ on $
(\RR^d\times\RR^d)$ with marginal $\eta$ w.r.t. the second coordinate and respectively  $\mu$ w.r.t. the first coordinate. 
Let  $\Ka_{n}$ be the Markov transitions on $\RR^d$ defined by the disintegration formulae
\begin{equation}\label{sinhorn-entropy-form-Sch}
\begin{array}{rclcrcl}
\Pa_{2n}(d(x,y))&=&\mu(dx)~\Ka_{2n}(x,dy)\in \Pi(\mu,\pi_{2n})& \mbox{\rm with}&
 \pi_{2n}&:=&\mu \Ka_{2n}\\
 &&&&&&\\
\Pa_{2n+1}(d(x,y))&=&\eta(dy)~\Ka_{2n+1}(y,dx)\in \Pi(\pi_{2n+1},\eta)&\mbox{\rm with}&
\pi_{2n+1}&:=&\eta \Ka_{2n+1}
\end{array}
\end{equation}
The collection of Sinkhorn transitions $\Ka_n$  are defined sequentially starting from a given reference measure 
$\Pa_0=P=\mu\times \Ka_0$. At every time step $n\geq 0$, given the distribution $\Pa_{2n}$ we choose the transition
$\Ka_{2n+1}$ as the $\Pa_{2n}$-conditional distribution of the first coordinate given the second. Given the distribution $\Pa_{2n+1}$ we choose the transition
 $\Ka_{2(n+1)}$ as the $\Pa_{2n+1}$-conditional distribution of the second coordinate given the first, and so on.
 
 Consider the disintegration of the Schr\"odinger bridge
 $$
 P_{\mu,\eta}=\mu\times \La_{\mu,\eta}
 \quad \mbox{\rm and}\quad
 \Pa_{2n}=\mu\times\Ka_{2n}
 $$
 Note that
$$
\Ha(  P_{\mu,\eta}~|~\Pa_{2n+1})=
\Ha(  P_{\mu,\eta}^{\flat}~|~\Pa_{2n+1}^{\flat})
$$
with the conjugate measures $( P^{\flat}_{\eta,\mu},\Pa^{\flat}_{2n+1})$ defined by the disintegration formulae
 \begin{eqnarray*}
 P^{\flat}_{\eta,\mu}(d(x,y))&:=&
  P_{\mu,\eta}(d(y,x))=\eta(dx)~ \La^{\flat}_{\eta,\mu}(x,dy)
 \\
\Pa^{\flat}_{2n+1}(d(x,y))&:=&\Pa_{2n+1}(d(y,x))=\eta(dx)~\Ka_{2n+1}(x,dy)
 \end{eqnarray*}
 \subsection{Some regularity conditions}
 We further assume that the reference transition $\Ka_0=K_W$  is the
linear Gaussian transition  associated with the potential function
\begin{equation}\label{def-W}
W(x,y)=-\log{g_{\tau}(y-(\alpha+\beta x))}
\end{equation}
 for a positive definite matrix $\tau>0$, some  $\alpha\in\RR^d$ and some invertible matrix $\beta\in\RR^{d\times d}$. In the above display, $g_{\tau}$ stands for the density of a centered Gaussian random variable with covariance matrix $\tau$. We underline that the linear Gaussian transition (\ref{def-W}) encapsulates all continuous time Gaussian models used in machine learning applications of Schr\" odinger bridges (see~\cite{adm-24,adm-25,dp} and references therein). Note that  $\nabla^2_2W(x,y)=\tau^{-1}$ so that $\Ka_{0}$  satisfies he $LS(\Vert \tau\Vert)$ and thus the $\TT_2(\Vert \tau\Vert)$-inequalities.

 Consider Boltzmann-Gibbs measures
 \begin{equation}\label{mu-eta}
\mu(dx):=e^{-U(x)}~\lambda(dx)\quad \mbox{\rm and}\quad
\eta(dx):=e^{-V(x)}~\lambda(dx)
\end{equation}
for some potential functions $(U,V)$. 
Assume there exists some $a_u,a_v>0$, $b_u,b_v\geq 0$ and $\delta_u,\delta_v\geq 0$ such that for any $x$ we have
 \begin{equation}\label{reg-UV}
  \nabla^2 U(x)\geq a_u~1_{\Vert x\Vert \geq  \delta_u}~I-b_u 
~1_{\Vert x\Vert< \delta_u}~I\quad \mbox{\rm and}\quad
  \nabla^2 V(x)\geq a_v~1_{\Vert x\Vert \geq  \delta_v}~I-b_v 
~1_{\Vert x\Vert< \delta_v}~I
\end{equation}
As shown in Section 7 in~\cite{dp}, the r.h.s. condition in (\ref{reg-UV}) ensures that $\La_{\mu,\eta}$ as well as Sinkhorn transitions $\Ka_{2n}$ with $n\geq 1$ and the Gibbs measure $\eta$ satisfy he $LS(\rho_v)$ and thus the $\TT_2(\rho_v)$-inequalities (\ref{K-log-sob}) with some parameter $\rho_v$ that depends on 
$(a_v,b_v,\delta_v)$. In addition, when $\delta_v=0$ we can choose $\rho_v=1/a_v$.

In the same vein, the l.h.s. condition in (\ref{reg-UV}) ensures that  $\La_{\eta,\mu}^{\flat}$ as well as
Sinkhorn transitions $\Ka_{2n+1}$ with $n\geq 1$ and the Gibbs measure $\mu$ satisfy the $LS(\rho_u)$ and thus the $\TT_2(\rho_u)$-inequalities (\ref{K-log-sob}) with some parameter $\rho_u$ that depends on 
$(a_u,b_u,\delta_u)$. In addition, when $\delta_u=0$ we can choose $\rho_u=1/a_u$.

Combining Theorem 7.9 with Remark 2.11 in~\cite{dp}  for any $n\geq 1$ we check following theorem.
\begin{theo}\label{theo-princ}
Consider the linear Gaussian transition  (\ref{def-W})  and assume the marginal potential functions $(U,V)$ satisfy (\ref{reg-UV}). In this situation, 
for any $n\geq 1$ we have 
the exponential decays
 \begin{equation}\label{def-eps}
 \Ha(P_{\mu,\eta}~|~\Pa_{n})\leq \left(1+ \varepsilon^{-1}
\right)^{-\lfloor n/2\rfloor}~ \Ha(P_{\mu,\eta}~|~\Pa_{0}) 
\quad \mbox{\rm with}\quad
\varepsilon:=\Vert \tau^{-1}\beta\Vert^2_2~\rho_u~\rho_v
\end{equation}
\end{theo}

Theorem 7.9~\cite{dp} provides exponential estimates for even and odd indices $n$ under slightly weaker regularity conditions.
For instance (\ref{def-eps}) holds for even indices $n=2k$ as soon as the  $V$ satisfies the r.h.s. condition in (\ref{reg-UV}) and $U$ satisfies the 
$\TT_2(\rho_u)$-inequality. 

When both  $U$ and $V$ are strongly log-concave (that is when (\ref{reg-UV}) holds with $\delta_v=0=\delta_v$), we can chose $\varepsilon:=\Vert \tau^{-1}\beta\Vert^2_2/(a_ua_v)$ in (\ref{def-eps}).
In this strongly log-concave  situation, Theorem 7.10~\cite{dp} provides sharper exponential convergence estimates.

To avoid unnecessary technical discussions on these different types of regularity conditions and possibly sharper convergence rates, the present article only discusses 
potential functions $(U,V)$ satisfying the curvature conditions (\ref{reg-UV}).

Last but not least, consider a Markov transition $\Ka$ satisfying the $LS(\rho)$ inequality for some $\rho>0$.
Choosing the reference measure $P=\mu\times \Ka$ in (\ref{ob-eq}) and
choosing $\eta=\pi\Ka$ for some probability measure $\pi$ we have $P_{\mu,\eta}=P_{\mu,\pi\Ka}=\mu\times \La_{\mu,\pi\Ka}$, as well as
$\mu\La_{\mu,\pi\Ka}=\pi\Ka$ and $P=P_{\mu,\mu\Ka}$. In this context,
 the estimate (\ref{eq-0}) takes the form
$$
\Da_2\left(\pi\Ka,\mu \Ka\right)^2\leq \int\mu(dx)~\Da_2\left(\delta_x\La_{\mu,\pi\Ka},\delta_x\Ka\right)^2
\leq 2\rho~\Ha\left(P_{\mu,\pi\Ka}~|~P_{\mu,\mu\Ka}\right)
$$
Assume there exists some parameter $\kappa>0$ such that for any $x,y\in\RR^d$ we have
$$
\Ja(\delta_x\Ka~|~\delta_y\Ka)\leq \kappa^2~\Vert x-y\Vert^2
$$
As shown in~\cite{dp}, all Sinkhorn transitions $\Ka=\Ka_n$ satisfies the above condition with the parameter $\kappa:=\Vert \tau^{-1}\beta\Vert_2$.
In this situation, under weaker conditions Theorem 1.2 in~\cite{dp} yields the following estimate
$$
2\Ha(P_{\mu,\pi\Ka}~|~P_{\mu,\mu\Ka})\leq \kappa^2\rho~ \Da_2(\pi,\mu)^2
$$
In addition, the bias and covariance entropic bounds estimates (\ref{eq-1}) and
(\ref{eq-2}) yields  
 \begin{eqnarray}
 \vertiii{m_{\Ka}-m_{\La_{\mu,\pi\Ka}}}_{2,\mu}&\leq &
 \kappa~\rho~ \Da_2(\pi,\mu)\label{eq-1-v2}\\
2^{-1} \vertiii{\sigma_{\Ka}-\sigma_{\La_{\mu,\pi\Ka}}}_{1,\mu}
&\leq& (\kappa\rho~\Da_2(\pi,\mu))^2+c_{\mu,\La}~\kappa~\rho~ \Da_2(\pi,\mu)\label{eq-2-v2}
 \end{eqnarray}

\subsection{Conditional mean and covariances}

To simplify notation we write $(m_{\mu,\eta},m_{2n})$ and $(\sigma_{\mu,\eta},\sigma_{2n})$  instead of  $(m_{\La_{\mu,\eta}},m_{\Ka_{2n}})$  and $(\sigma_{\La_{\mu,\eta}},\sigma_{\Ka_{2n}})$  the conditional mean and the conditional covariances associated with the bridge transitions
$\La_{\mu,\eta}
$ and Sinkhorn transitions $\Ka_{2n}$.

Using (\ref{eq-1}) and (\ref{eq-2}) we check the following bias and covariance estimates.
\begin{cor}\label{cor-even}
Consider the linear Gaussian transition  (\ref{def-W})  and assume the marginal potential functions $(U,V)$ satisfy (\ref{reg-UV}). In this situation, for any $n\geq 1$ we have
 \begin{eqnarray}
 \vertiii{m_{2n}-m_{\mu,\eta}}_{2,\mu}^2&\leq &
 2\rho_v~\Ha\left(P_{\mu,\eta}~|~\Pa_{2n}\right)\label{eq-1-e}\\
 \vertiii{\sigma_{2n}-\sigma_{\mu,\eta}}_{1,\mu}
&\leq& 4\rho_v~\Ha\left(P_{\mu,\eta}~|~\Pa_{2n}\right)+c_{\mu,\eta}~\left(8\rho_v~\Ha\left(P_{\mu,\eta}~|~\Pa_{2n}\right)\right)^{1/2}\label{eq-2-e}
 \end{eqnarray}
 with the parameter
$
c_{\mu,\eta}^2:=\int\mu(dx)\tr(\sigma_{\mu,\eta}(x))
$.
\end{cor}

As above, with a slight abuse of notation, we write $(m_{\eta,\mu},m_{2n+1})$ and $(\sigma_{\eta,\mu},\sigma_{2n+1})$  instead of  $(m_{\La^{\flat}_{\mu,\eta}},m_{\Ka_{2n+1}})$  and $(\sigma_{\La^{\flat}_{\eta,\mu}},\sigma_{\Ka_{2n+1}})$  the conditional mean and the conditional covariances associated with the transitions
$(\La^{\flat}_{\eta,\mu},\Ka_{2n+1})$. 

By (\ref{eq-1}) and (\ref{eq-2}) we have the following bias and covariance estimates.
\begin{cor}\label{cor-odd}
Consider the linear Gaussian transition  (\ref{def-W})  and assume the marginal potential functions $(U,V)$ satisfy (\ref{reg-UV}). In this situation, for any $n\geq 0$ we have
 \begin{eqnarray*}
 \vertiii{m_{2n+1}-m_{\eta,\mu}}_{2,\eta}^2&\leq &
 2\rho_u~\Ha\left(P_{\mu,\eta}~|~\Pa_{2n+1}\right)\label{eq-1-o}\\
 \vertiii{\sigma_{2n+1}-\sigma_{\eta,\mu}}_{1,\eta}
&\leq& 4\rho_u~\Ha\left(P_{\mu,\eta}~|~\Pa_{2n+1}\right)+c_{\eta,\mu}~\left(8\rho_u~\Ha\left(P_{\mu,\eta}~|~\Pa_{2n+1}\right)\right)^{1/2}\label{eq-2-o}
 \end{eqnarray*}
 with the parameter
$
c_{\eta,\mu}^2:=\int\eta(dx)\tr(\sigma_{\eta,\mu}(x))
$.
\end{cor}

Combining Theorem~\ref{theo-princ} with Corollary~\ref{cor-even} and Corollary~\ref{cor-odd}, several exponential decays estimates can be derived.
For instance, combining (\ref{def-eps}) and (\ref{eq-1-e}) we have
\begin{equation}\label{m2n-e}
 \vertiii{m_{2n}-m_{\mu,\eta}}_{2,\mu}^2\leq 
 2\rho_v~ \left(1+ \varepsilon^{-1}
\right)^{-n}~ \Ha(P_{\mu,\eta}~|~\Pa_{0}) 
\end{equation}
with $\varepsilon$ as in (\ref{def-eps}). In the same vein, combining (\ref{def-eps}) and (\ref{eq-2-e}) we have
\begin{equation}\label{s2n-e}
\begin{array}{l}
 \vertiii{\sigma_{2n}-\sigma_{\mu,\eta}}_{1,\mu}\\
 \\
\leq 4\rho_v~ \left(1+ \varepsilon^{-1}
\right)^{-n}~ \Ha(P_{\mu,\eta}~|~\Pa_{0}) +2c_{\mu,\eta}~\left(2\rho_v\right)^{1/2}~ \left(1+ \varepsilon^{-1}
\right)^{-n/2}~ \Ha(P_{\mu,\eta}~|~\Pa_{0})^{1/2}
\end{array}
\end{equation}
with $\varepsilon$ as in (\ref{def-eps}) and $c_{\mu,\eta}$ as in Corollary~\ref{cor-even}.
\subsection{Gradient and Hessian of Sinkhorn potentials}

Sinkhorn bridges $\Pa_n$ associated with the marginal measures (\ref{mu-eta}) and a reference transition of the form $K_W$ 
for some potential $W(x,y)$
are often expressed  in terms of a sequence of potential functions $(U_n,V_n)$ defined by the  integral recursions
\begin{eqnarray}
U_{2n}~&:=&U+\log{K_W(e^{-V_{2n}})}=U_{2n+1}\nonumber\\
V_{2n+1}&:=&V+\log{K_{W^{\flat}}(e^{-U_{2n+1}})}=V_{2(n+1)}
\label{prop-schp}
\end{eqnarray}
with $W^{\flat}(y,x):=W(x,y)$ and  the initial condition $V_0=0$.  In this notation, for any $n\geq 0$ we have
$$
\Pa_{n}(d(x,y))=e^{-U_n(x)}~e^{-W(x,y)}~e^{-V_n(y)}~\lambda(dx)\lambda(dy)
$$
In this situation, by Theorem 4.2 in~\cite{nutz} Schr\"odinger bridges  also take the form 
$$
P_{\mu,\eta}(d(x,y))=e^{-\UU(x)}~e^{-W(x,y)}~e^{-\VV(y)}~\lambda(dx)\lambda(dy)
$$
with the solution $(\UU,\VV)$ of the Schr\" odinger system
\begin{equation}
\UU=U+\log{K_W(e^{-\VV})}\quad \mbox{\rm and}\quad
\VV=V+\log{K_{W^{\flat}}(e^{-\UU})}
\label{prop-schp-fp}
\end{equation}
We assume that $(U,V)$ are twice differentiable. In this case, Sinkhorn potentials $(U_n,V_n)$ as well as  Schr\" odinger  potentials $(\UU,\VV)$ are also twice differentiable.

 We further assume that the reference transition $\Ka_0=K_W$  is  the linear Gaussian transition (\ref{def-W}). In this situation,
as shown in Example C.5 in~\cite{dp} we also have the gradient and Hessian formulae
  \begin{eqnarray}
\nabla U_{2n}(x)-\nabla \UU(x)&=&\cchi^{\prime}~(m_{2n}(x)-m_{\mu,\eta}(x))\nonumber\\
\nabla^2 U_{2n}(x)-\nabla^2\UU(x)&=&\cchi^{\prime}~(\sigma_{2n}(x)-\sigma_{\mu,\eta}(x))~\cchi\quad\mbox{\rm with}\quad
\cchi:=\tau^{-1}\beta\label{lem-u-mc}
 \end{eqnarray}
 This yields the following estimates.
 \begin{prop}\label{prop-even}
For any $n\geq 1$ we have
    \begin{eqnarray*}
\vertiii{\nabla U_{2n}-\nabla\UU}_{2,\mu}^2
&\leq &\Vert\cchi\Vert_2^2~\vertiii{m_{2n}-m_{\mu,\eta}}_{2,\mu}^2\\
\vertiii{\nabla^2 U_{2n}-\nabla^2\UU}_{1,\mu}&\leq& \Vert \cchi\Vert^2_2
~\vertiii{\sigma_{2n}-\sigma_{\mu,\eta}}_{1,\mu}
 \end{eqnarray*}
 \end{prop}
 Similarly, we have
  \begin{eqnarray}
\nabla V_{2n+1}(x)-\nabla \VV(x)&=&\cchi~(m_{2n+1}(x)-m_{\eta,\mu}(x))\nonumber\\
\nabla^2 V_{2n+1}(x)-\nabla^2 \VV(x)&=&\cchi~(\sigma_{2n+1}(x)-\sigma_{\eta,\mu}(x))~\cchi^{\prime}\label{lem-v-mc}
 \end{eqnarray}
 This yields the following estimates.
  \begin{prop}
For any $n\geq 0$ we have
    \begin{eqnarray*}
\vertiii{\nabla V_{2n+1}-\nabla\VV}_{2,\eta}^2
&\leq &\Vert\cchi\Vert_2^2~\vertiii{m_{2n+1}-m_{\eta,\mu}}_{2,\eta}^2\\
\vertiii{\nabla^2 V_{2n+1}-\nabla^2\VV}_{1,\eta}&\leq& \Vert \cchi\Vert^2_2
~\vertiii{\sigma_{2n+1}-\sigma_{\eta,\mu}}_{1,\eta}
 \end{eqnarray*}
 \end{prop}
 Combining the above propositions with Theorem~\ref{theo-princ}, Corollary~\ref{cor-even} and Corollary~\ref{cor-odd}, several exponential decays estimates can be derived.
 
 For instance, combining Proposition~\ref{prop-even} with  (\ref{m2n-e}) we check that
\begin{equation}\label{nablaU}
 \vertiii{\nabla U_{2n}-\nabla\UU}_{2,\mu}^2
\leq 2~\Vert\cchi\Vert_2^2~\rho_v~ \left(1+ \varepsilon^{-1}
\right)^{-n}~ \Ha(P_{\mu,\eta}~|~\Pa_{0}) 
\end{equation}
with $\varepsilon$ as in (\ref{def-eps}).
In the same vein, combining Proposition~\ref{prop-even} with  (\ref{s2n-e}) we check that
\begin{equation}\label{nabla2U}
\begin{array}{l}
\vertiii{\nabla^2 U_{2n}-\nabla^2\UU}_{1,\mu}\\
\\
\leq 4\Vert \cchi\Vert^2_2~
\rho_v~ \left(1+ \varepsilon^{-1}
\right)^{-n}~ \Ha(P_{\mu,\eta}~|~\Pa_{0}) +2c_{\mu,\eta}~\Vert \cchi\Vert^2_2~\left(2\rho_v\right)^{1/2}~ \left(1+ \varepsilon^{-1}
\right)^{-n/2}~ \Ha(P_{\mu,\eta}~|~\Pa_{0})^{1/2}
\end{array}
\end{equation}
with $\varepsilon$ as in (\ref{def-eps}) and $c_{\mu,\eta}$ as in Corollary~\ref{cor-even}.
Using the conditional covariance formulation  of the gradient and the Hessian of the bridges (\ref{lem-u-mc}), the estimates (\ref{nablaU}) and (\ref{nabla2U}) are direct consequences of the entropy continuity bounds presented in 
Theorem~\ref{theo-1} and the exponential entropy decays stated in 
Theorem~\ref{theo-princ}. 

Using a different approach, similar estimates in the case $(\alpha,\beta)=(0,I)$ and $\tau=tI$ in the asymptotics $t\downarrow 0$
are presented in the recent article~\cite{gt}. This kind of asymptotics is particularly useful as the conditional bridge expectations
converge to the Brenier map in the limit $t\downarrow 0$.

The proof of Theorem 1.2 in~\cite{gt} relies on the stochastic
control representation of entropic plans as laws of solutions to time-inhomogeneous stochastic differential equations and back propagated entropic potential along
Hamilton Jacobi equations
 in the spirit of~\cite{conforti-ptrf}. The central idea is to exploit the propagation of semiconcavity along Hamilton Jacobi equations to
obtain a quantitative stability result for primal solutions.
 This refined analysis allows to handle less regular models such as weakly log-concave marginals and situations where 
 the gradient/Hessian of the single potentials are not defined.

\subsection*{Acknowledgements}
The author warmly thanks Laurent Miclo for many useful discussions on log-Sobolev
and quadratic transportation cost
inequalities. He is also grateful to Giovanni Conforti for valuable comments on the regularity properties of Schr\"odinger bridges, including 
the stochastic
control representation of entropic plans and  backward propagations techniques of entropic potential along
Hamilton Jacobi equations.

\end{document}